\documentclass [a4paper,10pt]{article}
\usepackage [utf8] {inputenc}
\usepackage {amsmath}
\usepackage {amsthm}
\usepackage {amssymb}
\usepackage {amsfonts}
\newtheorem {lemma}{Lemma}
\newtheorem {theorem}{Theorem}

\title {The $L^1$ mean of the exponential sum of $d(n)$}
\author {Tomos Parry}
\date {}
\begin {document}
\maketitle 
\begin {abstract}
We use the classical circle method to give a relatively simple proof that
\[ \int _0^1\left |\sum _{n\leq x}d(n)e(n\alpha )\right |d\alpha \asymp \sqrt x.\]
\end {abstract}
\begin {center}
\section* {Introduction}
\end {center}
For a sequence $a_n$ (such as $\Lambda (n)$ or $d_k(n)$) let $S(\alpha)$ be its exponential sum 
\[ \sum_{n\leq x}a_ne(n\alpha )\]
and consider the $L^1$ norm
\[ \int _0^1|S(\alpha )|d\alpha .\]
We'd expect this norm to give us some insight about the sequence $a_n$ itself, the most famous result probably being Vaughan's result \cite {vaughanl1primes} that for the (weighted) primes the norm is $\gg \sqrt x$; the same argument shows the norm to be $\gg \sqrt x$ for the divisor functions too. On the other hand the Cauchy-Schwarz inequality shows the norm to be 
\begin {eqnarray*}
&\ll &\sqrt {x\log x}\hspace {10mm}\text { for the primes}
\\ &\ll &\sqrt {x(\log x)^{k^2-1}}\hspace {10mm}\text { for the $k$-fold divisor function}
\end {eqnarray*}
so the question is whether these $\log $ powers are really there. In this note we'll be interested in this question for the divisor function $d(n)$.
\\
\\ That the norm is $\ll \sqrt x$ was claimed in \cite {anghywir}, but the proof was shown to be incorrect by Goldston and Pandey \cite {goldstonpandey}, who also gave the current best bound in the literature - they proved the norm is $\ll \sqrt {x}\log x$. We prove that it is $\ll \sqrt x$.
\\
\\ Parts of this work was started when the author was working for the Indraprastha Institute of Information Technology, Delhi, and we thank Prof.'s Sneha Chaubey and Debika Benerjee for discussions around the problem.
\\ \begin {theorem}\label {t1}
\[ \text {If }\hspace {5mm}S(\alpha )=\sum _{n\leq x}d(n)e(n\alpha )\hspace {10mm}\text {then}\hspace {5mm}\int _0^1\left |S(\alpha )\right |d\alpha \asymp \sqrt x.\]
\end {theorem}
The lower bound is already given in \cite {goldstonpandey} so we only prove the upper bound. We prove it subejct to three lemmas which we prove in the separate sections after. Write
\[ \zeta (s)=\frac {1}{s-1}+a_1+a_2(s-1)+\cdot \cdot \cdot \hspace {10mm}.\]
First we need the explicit coefficients in evaluating $\sum _{n\leq x}d(n)^2$.
\begin {lemma}\label {diagonal}
Write $\mathcal F(s)=1/\zeta (2s)$ and let
\[ c_{3,0}=\frac {1}{6}\hspace {5mm}c_{2,1} =\frac {1}{2}\hspace {5mm}c_{2,0}=2a_1-\frac {1}{2}\hspace {5mm}c_{1,2}=\frac {1}{2}\hspace {5mm}c_{1,1}=4a_1-1\hspace {5mm}c_{1,0}=4a_2+6a_1^2-4a_1+1.\]
Then
\[ \sum _{n\leq x}d(n)^2=x\sum _{J+K\leq 3}(\log x)^J\mathcal F^{(K)}(1)c_{J,K}+\mathcal O\left (x^{1/2+\epsilon }\right ).\]
\end {lemma}
We then need the explicit coefficients for a sum coming from an exponential sum approximation. Let $\gamma =x^{1/\Delta }$.
\begin {lemma}\label {ail}
Write $\mathcal G(s)=1/\zeta (s+1)$ and let
\begin {eqnarray*}
&&d_{3,0}=\frac {4}{3\Delta ^3}-\frac {2}{\Delta ^2}+\frac {1}{\Delta }\hspace {5mm}d_{2,1}=1\hspace {5mm}d_{2,0}=a_1+\frac {2(2a_1-1)}{\Delta }\left (1-\frac {1}{\Delta }\right )
\\ &&d_{1,2}=2\hspace {5mm}d_{1,1}=8a_1-2\hspace {5mm}d_{1,0}=4a_2+4a_1^2-2a_1+\frac {4a_1^2-4a_1+2}{\Delta }.
\end {eqnarray*}
Then 
\[ \sum _{q\leq \gamma }\frac {\phi (q)}{q^2}\int _1^xg_q(t)^2dt=x\sum _{J+K\leq 3}(\log x)^J\mathcal G^{(K)}(1)d_{J,K}+\mathcal O\left (\frac {x}{\sqrt \gamma }\right ).\]
\end {lemma}
And we'll need to deal with a cross-term, for which we'll need an equidistribution result for $d(n)$ on average in arithmetic progressions of moduli up to $\sqrt x$. This is certainly possible since even in individual progressions $d(n)$ is evenly distributed up to $x^{2/3}$, but we have to avoid $x^\epsilon $ factors. 
\begin {lemma}\label {3}
Let 
\begin {eqnarray}\label {m}
\mathcal M_x(q,a)&:=&\frac {x}{q}\sum _{r|q}\frac {c_r(a)}{r}\Big (\log (x/r^2)+2\gamma -1\Big )
\\ E_x(q,a)&:=&\sum _{n\leq x\atop {n\equiv a\text {mod }(q)}}d(n)-\mathcal M_x(q,a).\notag 
\end {eqnarray}
Then for $q\leq \sqrt x$
\begin {eqnarray*}
\sum _{a=1}^q|E_x(q,a)|^2&\ll &q\sqrt x.
\end {eqnarray*}
\end {lemma}
Now we can do the upper bound in Theorem \ref {t1}. Let $\gamma >0$ and for $(a,q)=1$ define
\begin {eqnarray*}
\mathfrak F(a/q)&=&\left (\frac {a}{q}-\frac {1}{q\overline a},\frac {a}{q}+\frac {1}{q\overline {(-a)}}\right )\hspace {10mm}\overline a\text { chosen in }(\gamma -q,\gamma ]
\end {eqnarray*}
so that for any continuous $f:\mathbb R\rightarrow \mathbb C$ of period 1 
\begin {equation}\label {farey}
\int _0^1f(\alpha )d\alpha =\sum _{q\leq \gamma }\sideset {}{'}\sum _{a=1}^q\int _{\mathfrak F(a/q)}f(\alpha )d\alpha \hspace {10mm}\text {length of }\mathfrak F(a/q)\asymp \frac {1}{q\gamma }.
\end {equation}
This is the \emph {Farey dissection} of the unit interval of order $\gamma$ and for its discussion see, for example, Section 3.8 of \cite {hardywright}. For $\alpha \in \mathfrak F(a/q)$ write always $\alpha =a/q+\beta $ with $(a,q)=1$. Let
\begin {eqnarray}\label {coffi}
f_q(t)&=&\log (t/q^2)+2\gamma -1\hspace {7mm}g_q(t)=\frac {d}{dt}\left \{ tf_q(t)\right \} \hspace {7mm}I_q(\beta )=\int _1^xe(t\beta )g_q(t)\notag 
\\ S_t(\alpha )&:=&\sum _{n\leq t}d(n)e(\alpha n)\hspace {10mm}S_t^*(\alpha )=\frac {I_q(\beta )}{q}\hspace {7mm}\Delta _t(\alpha )=S_t(\alpha )-S_t^*(\alpha )
\end {eqnarray}
so that 
\begin {eqnarray}\label {bora}
\int _0^1\left |\Delta (\alpha )\right |^2d\alpha &=&\underbrace {\int _0^1|S(\alpha )|^2d\alpha }_{=:\mathcal M_1}-\underbrace {\int _0^1|S^*(\alpha )|^2d\alpha }_{=:\mathcal M_2}+\mathcal O\left (\int _0^1\left |S^*(\alpha )\Delta (\alpha )\right |d\alpha \right );
\end {eqnarray}
note also that integration by parts gives
\begin {eqnarray}\label {afal}
S_x^*(\alpha )&=&\frac {1}{q}\left (xf_x(q)-2\pi i\int _1^xe(t\beta )tf_t(q)dt\right )+\mathcal O(1)
\end {eqnarray}
as well as 
\begin {eqnarray}
I_q(\beta )&=&\mathcal O\left (\int _1^{q^2}|g_q(t)|dt\right )+\int _{q^2}^xe(\beta t)g_q(t)dt
\hspace {2mm}\ll \hspace {2mm}q^2+\frac {\log (x/q^2)}{\beta }\hspace {10mm}\text {for all $\beta \not =0$}\label {i2}
\\ &\ll &q^2+\frac {\log (x/q^2)x}{1+|\beta |x}\hspace {10mm}\text {for $\alpha \in \mathfrak F(a/q)$ by \eqref {farey}}\label {i}
\end {eqnarray}
which with \eqref {farey} in particular says
\begin {eqnarray}\label {ola}
\sum _{q\leq \gamma }\sideset {}{'}\sum _{a=1}^q\int _{\mathfrak F(a/q)}|S_x^*(\alpha )|d\alpha &\ll &\sum _{q\leq \gamma }\log (x/q^2)\hspace {2mm}\ll \hspace {2mm}\gamma \log (x/\gamma ^2).
\end {eqnarray}
From \eqref {coffi}, assuming $(a,q)=1$,
\begin {eqnarray}\label {hwn}
S_x^*(a/q)&=&\frac {x}{q}f_x(q)\hspace {2mm}=:\hspace {2mm}F(q)
\end {eqnarray}
so that \eqref {m} and \eqref {coffi} say,
\begin {eqnarray*}
\mathcal M_x(q,a)&=&\frac {1}{q}\sum _{r|q}\frac {c_r(a)}{r}F(r)
\end {eqnarray*}
so, not assuming $(b,q)=1$,
\begin {eqnarray*}
\sum _{a=1}^q\mathcal M_x(q,a)e\left (\frac {ab}{q}\right )&=&F\left (\frac {q}{(q,b)}\right )\hspace {2mm}=\hspace {2mm}S_x^*(b/q)
\end {eqnarray*}
so from \eqref {coffi} and with $E_x(q,a)$ as in Lemma \ref {3}
\begin {eqnarray*}
E_x(q,a)=\frac {1}{q}\sum _{b=1}^qe\left (-\frac {ab}{q}\right )\left (\sum _{n\leq x}d(n)e\left (\frac {nb}{q}\right )-\sum _{N=1}^q\mathcal M_x(q,N)e\left (\frac {Nb}{q}\right )\right )=\frac {1}{q}\sum _{b=1}^qe\left (-\frac {ab}{q}\right )\Delta _x(b/q)
\end {eqnarray*}
and therefore 
\begin {eqnarray}\label {11}
\sum _{b=1}^q\left |\Delta _x(b/q)\right |^2&=&q\sum _{a=1}^q|E_{x}(q,a)|^2\hspace {2mm}\ll q^2\sqrt x\hspace {10mm}\text { for }q\leq \sqrt x\hspace {2mm}
\end {eqnarray}
by Lemma \ref {3}. By partial summation and \eqref {coffi} 
\begin {eqnarray*}
S_x(\alpha )&=&e(x\beta )S_x(a/q)-2\pi i\beta \int _1^xe(\beta t)S_t(a/q)dt
\end {eqnarray*}
so by \eqref {coffi} and \eqref {afal}
\begin {eqnarray}\label {22}
\Delta _x(\alpha )&=&e(x\beta )\left (S_x(a/q)-\frac {x}{q}f_q(x)\right )-2\pi i\beta \int _1^xe(t\beta )\left (S_t(a/q)-\frac {tf_t(q)}{q}\right )dt+\mathcal O(1)\notag 
\\ &\ll &\left (1+|\beta |x\right )\max _{t\leq x}|\Delta _t(a/q)|\hspace {10mm}\text { (by \eqref {hwn} and \eqref {coffi})}\notag 
\end {eqnarray}
so with \eqref {11} we get
\begin {eqnarray*}
\sum _{a=1}^q|\Delta _x(\alpha )|&\ll &(1+|\beta |x)q^{3/2}x^{1/4}
\end {eqnarray*}
so from \eqref {i}
\begin {eqnarray*}
\sum _{a=1}^q\int _{\pm 1/q\gamma }|I_q(\beta )\Delta _x(\alpha )|d\beta &\ll &\log \left (\frac {x}{q\gamma }\right )\frac {q^{1/2}x^{5/4}}{\gamma } 
\end {eqnarray*}
and with \eqref {farey}, \eqref {coffi} we conclude the error term in \eqref {bora} is
\begin {eqnarray*}
\sum _{q\leq \gamma }\frac {1}{q}\sideset {}{'}\sum _{a=1}^q\int _{\mathfrak F(a/q)}\left |I_q^*(\alpha )\Delta _x(\alpha )\right |d\alpha &\ll &\frac {x^{5/4}}{\gamma }\sum _{q\leq \gamma }\frac {\log (x/q\gamma )}{\sqrt q}\hspace {2mm}\ll \hspace {2mm}\frac {x^{5/4}}{\sqrt \gamma }\log \left (\frac {x}{\gamma ^2}\right ).
\end {eqnarray*}
As for the main terms in \eqref {bora}, first let
\begin {eqnarray*}
L(Z)&=&\int _{1,1-Z}^{x,x-Z}g_q(t)g_q(t-Z)dt\hspace {7mm}\hat L(\beta )=\int _{\pm \infty }e(-\beta Z)L(Z)dZ
\\ c_{J,K}&=&\text { as in Lemma \ref {diagonal}}\hspace {7mm}d_{J,K}=\text { as in Lemma \ref {ail}}
\end {eqnarray*}
so that $|I_q(\beta )|^2=\hat L(\beta )$ and therefore
\begin {eqnarray}\label {fourier}
\int _{\pm \infty }|I_q(\beta )|^2d\beta =L(0).
\end {eqnarray}
By \eqref {farey} and \eqref {i2}
\[ \left (\int _{\pm \infty }-\int _{\mathfrak F(a/q)}\right )|I_q(\beta )|^2d\beta \ll q\gamma \left (\log \left (\frac {x}{q^2}\right )\right )^2\]
which with \eqref {farey} and \eqref {fourier} gives (referring to the definitions in \eqref {bora} and \eqref {coffi})
\begin {eqnarray*}
\mathcal M_2&=&\sum _{q\leq \gamma }\frac {1}{q^2}\sideset {}{'}\sum _{a=1}^q\int _{\mathfrak F(a/q)}|I_q(\beta )|^2d\alpha 
\\ &=&\sum _{q\leq \gamma }\frac {\phi (q)}{q^2}\int _{\pm \infty }|I_q(\beta )|^2d\beta +\mathcal O\left (\gamma \sum _{q\leq \gamma }\left (\log (x/q^2)\right )^2\right )
\\  &=&\sum _{J+K\leq 3}(\log x)^J\mathcal G^{(K)}(1)d_{J,K}+\mathcal O\left (\frac {x}{\sqrt \gamma }+\gamma ^2\log (x/\gamma ^2)\right )
\end {eqnarray*}
by Lemma \ref {ail}, whilst 
\begin {eqnarray*}
\mathcal M_1&=&x\sum _{J+K\leq 3}(\log x)^J\mathcal F^{(K)}(1)c_{J,K}+\mathcal O\left (x^{1/2+\epsilon }\right )
\end {eqnarray*}
by Lemma \ref {diagonal}. All in all
\begin {eqnarray*}
\int _0^1\left |\Delta (\alpha )\right |^2d\alpha &=&x\sum _{J+K\leq 3}(\log x)^J\left (\mathcal F^{(K)}(1)c_{J,K}-\mathcal G^{(K)}(1)d_{J,K}\right )+\mathcal O\left (\left (\frac {x^{5/4}}{\sqrt \gamma }+\gamma ^2\right )\right )
\\ &=&x\sum _{J+K\leq 3}(\log x)^J\mathcal F^{(K)}(1)\left (c_{J,K}-\frac {d_{J,K}}{2^K}\right )+\mathcal O\left (\left (\frac {x^{5/4}}{\sqrt \gamma }+\gamma ^2\right )\log (x/\gamma ^2)\right ).
\end {eqnarray*}
Now it seems something important happens - if we set $\Delta =2$ then the $c_{J,K}$ and $d_{J,K}/2^K$ coefficients match up. So actually
\[ \int _0^1|\Delta (\alpha )|^2d\alpha \ll x\]
and with \eqref {ola} we conclude 
\begin {eqnarray*}
\int _0^1|S(\alpha )|d\alpha &=&\int _0^1|S^*(\alpha )|d\alpha +\mathcal O\left (\left (\int _0^1|\Delta (\alpha )|^2d\alpha \right )^{1/2}\right )\ll \sqrt x
\end {eqnarray*}
and we're done.
\begin {center}
\section* {Proof of Lemma \ref {diagonal}}
\end {center}
Write
\[ \zeta (s)=\frac {1}{s-1}+a_1+a_2(s-1)+\cdot \cdot \cdot \hspace {10mm}l_N(n)=\sum _{k_1+\cdot \cdot \cdot +k_{N}=n}a_{k_1}\cdot \cdot \cdot a_{k_N}\]
so that, for any $A,N,X\geq 0$, 
\begin {eqnarray}\label {dechra}
\lim _{s\rightarrow 1}\left (\frac {d}{ds}\right )^A\Big \{ (s-1)^N\zeta (s)^N\Big \} &=&A!l_NA
\\ \lim _{s\rightarrow 1}\left (\frac {d}{ds}\right )^{X}\left \{ (s-1)^{A+1}\zeta ^{A}(s)\right \} &=&X!(X-1)\cdot \cdot \cdot (X-A)a_X\left \{ \begin {array}{ll}1&\text { if }X\geq A+1\\ 0&\text { otherwise }\end {array}\right .
\\ \lim _{s\rightarrow 1}\left (\frac {d}{ds}\right )^{X}\left \{ (s-1)^{A+1}\zeta ^{A}(s)\right \} &=&(-1)^AA!\hspace {10mm}(X=0)
\end {eqnarray}
so
\begin {eqnarray*}
\sum _{n\leq x}d(n)^2&=&Res_{s=1}\left \{ \frac {\zeta (s)^4\mathcal F(s)x^s}{s}\right \} 
=x\sum _{B+D\leq 3}\frac {(\log x)^B}{B!D!}\mathcal F^{(D)}(-1)^{3-B-D}\sum _{0\leq n\leq 3-B-D}(-1)^nl_4(n)
\end {eqnarray*}
and the claim follows.
\begin {center}
{\section* {Proof of Lemma \ref {ail}}}
\end {center}
First we prove
\newtheorem* {lemmaa}{Lemma \ref {ail}A}
\begin {lemmaa}\label {ok}
For any $\alpha _{0,0},\alpha _{0,1},\alpha _{1,0}\in \mathbb C$ let
\begin {align*}
\beta _{Q,X}&=\left \{ \begin {array}{ll}
\alpha _{0,0}+\alpha _{0,1}&\text { if }X+Q=0
\\ \alpha _{Q,X}&\text { if }X+Q=1.
\end {array}\right \} &g_q(t)&=\sum _{Q+X\leq 1}\beta _{Q,X}(\log q)^Q(\log t)^X
\\ S(\mathcal Q^*,\mathcal X)&=\sum _{\mathbf Q+\mathbf X\leq 1\atop {Q+Q'+1=\mathcal Q^*\atop {X+X'=\mathcal X}}}\left (\alpha _{\mathbf Q,\mathbf X}+\underbrace {2\alpha _{0,1}\alpha _{Q,X}}_{Q'=X'=0}+\underbrace {\alpha _{\mathbf 0,\mathbf 1}}_{\mathbf Q=\mathbf X=0}\right )&\mu _{\mathcal Q^*,\mathcal X}&=-\frac {a_{\mathcal Q^*-K}(-1)^{J+\mathcal Q^*+\mathcal X}(\mathcal Q^*-1)!\mathcal X!}{J!K!}\cdot
\\ t(J,K)&=\sum _{1\leq \mathcal Q^*+\mathcal X\leq 3\atop {\mathcal X\geq J\atop {\mathcal Q^*\geq K}}}\mu _{\mathcal Q^*,\mathcal X}S(\mathcal Q^*,\mathcal X)&\gamma _{\mathcal Q^*,\mathcal X}&=\mu _{\mathcal Q^*,\mathcal X}\left \{ \begin {array}{ll}
1&\mathcal Q^*=K
\\ 0&\mathcal Q^*>K
\end {array}\right \} 
\\ t(J)&=\sum _{J\leq \mathcal Q^*+\mathcal X\leq 3\atop {\mathcal Q^*\leq J}}\gamma _{\mathcal Q^*,\mathcal X}^*S(\mathcal Q^*,\mathcal X)&\gamma ^*_{\mathcal Q^*,\mathcal X}&=\frac {(-1)^{J+\mathcal Q^*+\mathcal X}\mathcal X!}{(J-\mathcal Q^*)!\mathcal Q^*\Delta ^{\mathcal Q^*}}.
\end {align*}
Then for $\gamma =x^{1/\Delta }$
\begin {eqnarray*}
\frac {1}{x}\sum _{q\leq \gamma }\frac {\phi (q)}{q^2}\int _1^xg_q(t)^2dt&=&\sum _{J+K\leq 3}(\log x)^J\mathcal G^{(K)}(1)t(J,K)+\mathcal G^{}(1)\sum _{J\leq 3}(\log x)^Jt(J)+\mathcal O\left (\frac {1}{\sqrt \gamma }\right ).
\end {eqnarray*}
\end {lemmaa}
\begin {proof}[Proof of Lemma 2A]
Since
\[ \int _1^x(\log t)^n=x\sum _{r=0}^n\frac {n!(-1)^r(\log x)^{n-r}}{(n-r)!}+\mathcal O(1)\hspace {10mm}(n\geq 0)\]
we have
\begin {eqnarray*}
\int _1^xg_q(t)^2dt&=&x\sum _{r+Q+Q'\leq 2}\alpha _{r,\mathbf Q}(\log x)^r(\log q)^{Q+Q'}+\mathcal O(1)\hspace {10mm}\alpha _{r,\mathbf Q}:=\sum _{\mathbf X\atop {\mathbf Q+\mathbf X\leq 1\atop {X+X'\geq r}}}\frac {(-1)^{X+X'-r}(X+X')!\beta _{\mathbf Q,\mathbf X}}{r!}
\end {eqnarray*}
so as (writing $\mathcal Q=Q+Q'$)
\[ \sum _{q\leq \gamma }\frac {\phi (q)(-\log q)^{\mathcal Q}}{q^2}=Res_{s=1}\left \{ \frac {\mathcal A^{(\mathcal Q)}(s)\gamma ^{s-1}}{s-1}\right \} +\mathcal O\left (\frac {1}{\sqrt \gamma }\right )\hspace {10mm}\mathcal A(s)=\zeta (s)\mathcal G(s)\]
we get up to an acceptable error
\begin {eqnarray*}
\sum _{q\leq \gamma }\frac {\phi (q)}{q^2}\int _1^xg_q(t)^2dt&=&x\sum _{r+Q+Q'\leq 2}(-1)^\mathcal Q(\log x)^r\underbrace {\alpha _{r,\mathbf Q}Res_{s=1}\left \{ \frac {\mathcal A^{(\mathcal Q)}(s)\gamma ^{s-1}}{s-1}\right \} }_{=:(\star )}.
\end {eqnarray*}
From \eqref {dechra}
\begin {eqnarray*}
&&\lim _{s\rightarrow 1}\left \{ \left (\frac {d}{ds}\right )^X(s-1)^{A+1}\zeta ^{A}(s)\left (\frac {d}{ds}\right )^Y\mathcal G^{B}(s)\left (\frac {d}{ds}\right )^Z\gamma ^{s-1}\right \} 
\\ &&\hspace {10mm}=\hspace {4mm}\underbrace {\mathcal G^{(B+Y)}(1)}_{=:\mathcal G^{B+Y}}(\log \gamma )^Z\left \{ \begin {array}{ll}X!(X-1)\cdot \cdot \cdot (X-A)a_X&X\geq A+1\hspace {2mm}X\not =0\\ 0&X\leq A\hspace {2mm}X\not =0\\ (-1)^AA!&X=0\end {array}\right .
\end {eqnarray*}
so $(\star )$ is
\begin {eqnarray*}
&=&\alpha _{r,\mathbf Q}\sum _{A+B=\mathcal Q}{{\mathcal Q}\choose {A,B}}A!\left ((-1)^A\sum _{Y+Z=A+1}\frac {\mathcal G^{B+Y}(\log x)^Z}{Y!Z!\Delta ^Z}+a_{A+1}\mathcal G^B\right )
\\ &=:&\sum _{A+B=\mathcal Q}\sum _{Y+Z=A+1}\mathcal G^{B+Y}X_{\mathbf Q,Y}^{r,A,B,Z}(\log x)^Z+\sum _{A+B=\mathcal Q}\mathcal G^BY_{r,B,\mathbf Q}^{A}
\end {eqnarray*}
so
\begin {eqnarray*}
&&\sum _{r+Q+Q'\leq 2}(-1)^\mathcal Q(\log x)^r\alpha _{r,\mathbf Q}Res_{s=1}\left \{ \frac {\mathcal A^\mathcal Q(s)\gamma ^{s-1}}{s-1}\right \} 
\\ &&\hspace {10mm}=\hspace {4mm}\sum _{J+K\leq 3}(\log x)^J\mathcal G^K\sum _{\mathcal Q,Y\atop {0,K-\mathcal Q\leq Y\leq K\leq \mathcal Q^*\leq J+K}}(-1)^\mathcal QX_{\mathbf Q,Y}^{J+K-\mathcal Q^*,\mathcal Q-K+Y,K-Y,\mathcal Q^*-K}
\\ &&\hspace {20mm}+\hspace {4mm}\sum _{J+K\leq 2}(\log x)^J\mathcal G^K\sum _{K\leq \mathcal Q\leq 2-J}(-1)^\mathcal QY_{J,K,\mathbf Q}^{\mathcal Q-K}
\\ &&\hspace {10mm}=:\hspace {4mm}\sum _{J+K\leq 3}(\log x)^J\mathcal G^KZ(J,K).
\end {eqnarray*}
The $\mathcal Q,Y$ sum is
\begin {eqnarray*}
&&\sum _{\mathbf Q+\mathbf X\leq 1\atop {K\leq \mathcal Q^*\leq J+K\leq \mathcal Q^*+\mathcal X}}\frac {(-1)^{J+\mathcal Q^*+\mathcal X}\mathcal Q!\mathcal X!}{(J+K-\mathcal Q^*)!(\mathcal Q^*-K)!\Delta ^{\mathcal Q^*-K}}\beta _{\mathbf Q,\mathbf X}\underbrace {\sum _{0,K-\mathcal Q\leq Y\leq K}\frac {(-1)^Y}{(K-Y)!Y!}}_{=0\text { unless }\mathcal Q^*=K\text { or }K=0}
\\ &&\hspace {20mm}=\hspace {4mm}\sum _{\mathbf Q+\mathbf X\leq 1\atop {K\leq \mathcal Q^*\leq J+K\leq \mathcal Q^*+\mathcal X}}\left \{ \begin {array}{ll}\gamma _{\mathcal Q^*,\mathcal X}&K\not =0\\ \gamma _{\mathcal Q^*,\mathcal X}^*&K=0\end {array}\right .
\end {eqnarray*}
and the $\mathcal Q$ sum is
\begin {eqnarray*}
\sum _{\mathbf Q+\mathbf X\leq 1\atop {\mathcal X\geq J\atop {K\leq \mathcal Q\leq 2-J}}}\frac {a_{\mathcal Q^*-K}(-1)^{J+\mathcal Q+\mathcal X}\mathcal Q!\mathcal X!\beta _{\mathbf Q,\mathbf X}}{J!K!}=\sum _{\mathbf Q+\mathbf X\leq 1\atop {\mathcal X\geq J\atop {\mathcal Q\geq K}}}\mu _{\mathcal Q^*,\mathcal X}\beta _{\mathbf Q,\mathbf X}
\end {eqnarray*}
so 
\begin {eqnarray*}
Z(J,K)&=&\sum _{\mathbf Q+\mathbf X\leq 1\atop {K\leq \mathcal Q^*\leq J+K\leq \mathcal Q^*+\mathcal X}}\beta _{\mathbf Q,\mathbf X}\left \{ \begin {array}{ll}\gamma _{\mathcal Q^*,\mathcal X}&K\not =0\\ \gamma _{\mathcal Q^*,\mathcal X}^*&K=0\end {array}\right \} +\sum _{\mathbf Q+\mathbf X\leq 1\atop {\mathcal X\geq J\atop {\mathcal Q^*>K}}}\beta _{\mathbf Q,\mathbf X}\mu _{\mathcal Q^*,\mathcal X}
\\ &=&\underbrace {\sum _{\mathbf Q+\mathbf X\leq 1\atop {\mathcal Q^*\leq J\leq \mathcal Q^*+\mathcal X}}\beta _{\mathbf Q,\mathbf X}\gamma _{\mathcal Q^*,\mathcal X}^*}_{K=0}+\sum _{\mathbf Q+\mathbf X\leq 1\atop {\mathcal X\geq J\atop {\mathcal Q^*\geq K}}}\beta _{\mathbf Q,\mathbf X}\mu _{\mathcal Q^*,\mathcal X}.
\end {eqnarray*}
Then
\begin {eqnarray*}
\sum _{\mathbf Q+\mathbf X\leq 1\atop {\mathcal Q^*\leq J\leq \mathcal Q^*+\mathcal X}}\beta _{\mathbf Q,\mathbf X}\gamma _{\mathcal Q^*,\mathcal X}^*&=&\sum _{\mathbf Q+\mathbf X\leq 1\atop {\mathcal Q^*\leq J\leq \mathcal Q^*+\mathcal X}}\alpha _{\mathbf Q,\mathbf X}\gamma _{\mathcal Q^*,\mathcal X}^*+2\alpha _{0,1}\sum _{Q+X\leq 1\atop {Q+1\leq J\leq Q+X+1}}\alpha _{Q,X}\gamma _{Q+1,X}^*+\underbrace {\alpha _{0,1}^2\gamma _{1,0}^*}_{1\leq J\leq 1}
\\ &=&\sum _{1\leq \mathcal Q^*+\mathcal X\leq 3\atop {\mathcal Q^*\leq J\leq \mathcal Q^*+\mathcal X}}\gamma _{\mathcal Q^*,\mathcal X}^*S(\mathcal Q^*,\mathcal X)
\\ &=&t(J)
\end {eqnarray*}
with a similar expression for $\mu _{\mathcal Q^*,\mathcal X},t(J,K)$ too, and the claim follows.
\end {proof}
To prove Lemma \ref {ail} we now need to calculate explicitly the $t(J,K)$'s and $t(J)$'s in Lemma 2A. There are quite a few values to calculate, but if we wait for a rainy Sunday we can draw up the tables on the next page:
\newpage
\begin{center}
\begin{tabular}{||c c c c||} 
 \hline
$(J,K)$ & $(\mathcal Q^*,\mathcal X)$ & $\mu_{\mathcal Q^*,\mathcal X}$ & $\gamma _{\mathcal Q^*,\mathcal X}^*$ \\ [0.5ex] 
 \hline\hline 
$(3,0)$ & $(3,0)$ & & $1/3\Delta ^3$ \\ 
 \hline
$(3,0)$ & $(2,1)$ & & $1/2\Delta ^2$ \\ 
 \hline
$(3,0)$ & $(1,2)$ & & $1/\Delta $ \\ 
 \hline
$(2,1)$ & $(1,2)$ & $1$ & \\ 
 \hline
$(2,0)$ & $(1,1)$ & & $1/\Delta $ \\ 
 \hline
$(2,0)$ & $(1,2)$ & $a_1$ & $-2/\Delta $\\ 
 \hline
$(2,0)$ & $(2,0)$ & & $1/2\Delta ^2$ \\ 
 \hline
$(2,0)$ & $(2,1)$ & & $-1/2\Delta ^2$ \\ 
 \hline
$(1,2)$ & $(2,1)$ & $-1/2$ & \\
 \hline
$(1,1)$ & $(1,1)$ & $1$ & \\
 \hline
$(1,1)$ & $(1,2)$ & $-2$ & \\
 \hline
$(1,1)$ & $(2,1)$ & $-a_1$ & \\
 \hline 
$(1,0)$ & $(1,0)$ &  & $1/\Delta $ \\
 \hline 
$(1,0)$ & $(1,1)$ & $a_1$ & $-1/\Delta $ \\
 \hline 
$(1,0)$ & $(1,2)$ & $-2a_1$ & $2/\Delta $ \\ 
 \hline 
$(1,0)$ & $(2,1)$ & $-a_2$ &  \\ 
 \hline
\end{tabular}
\hspace {10mm}
\begin{tabular}{||c c ||} 
 \hline
$(\mathcal Q^*,\mathcal X)$ & $S(\mathcal Q^*,\mathcal X)$ \\ [0.5ex] 
 \hline\hline
$(3,0)$ & $\alpha _{1,0}^2$  \\ 
 \hline
$(2,1)$ & $2\alpha _{1,0}\alpha _{0,1}$  \\ 
 \hline
$(2,0)$ & $2\alpha _{1,0}\left (\alpha _{0,0}+\alpha _{0,1}\right )$ \\ 
 \hline
$(1,2)$ & $\alpha _{0,1}^2$ \\
 \hline
$(1,1)$ & $2\alpha _{0,1}\left (\alpha _{0,0}+\alpha _{0,1}\right )$ \\
 \hline 
$(1,0)$ & $(\alpha _{0,0}+\alpha _{0,1})^2$ \\ 
 \hline
\end{tabular}
\end{center}
\begin{center}
\begin{tabular}{||c c c ||} 
 \hline
$(J,K)$ & $t(J,K)$ & $t(J)$ \\ [0.5ex] 
 \hline\hline
$(3,0)$ &$0$& $\gamma _{3,0}^*S(3,0)+\gamma _{2,1}^*S(2,1)+\gamma _{1,2}^*S(1,2)$  \\ 
 \hline
$(2,1)$ & $\mu _{1,2}S(1,2)$ & \\ 
 \hline
$(2,0)$ &$\mu _{1,2}S(1,2)$& $\gamma _{1,1}^*S(1,1)+\gamma _{1,2}^*S(1,2)+\gamma _{2,0}^*S(2,0)+\gamma _{2,1}^*S(2,1)$ \\ 
 \hline
$(1,2)$ & $\mu _{2,1}S(2,1)$ & \\
 \hline
$(1,1)$ & $\mu _{1,1}S(1,1)+\mu _{1,2}S(1,2)+\mu _{2,1}S(2,1)$ & \\
 \hline 
$(1,0)$ &$\mu _{1,1}S(1,1)+\mu _{1,2}S(1,2)+\mu _{2,1}S(2,1)$& $\gamma _{1,0}^*S(1,0)+\gamma _{1,1}^*S(1,1)+\gamma _{1,2}^*S(1,2)$\\
 \hline
\end{tabular}
\end{center}
Recall $\alpha _{0,0}=2a_1-1$ $\alpha _{0,1}=1$ $\alpha _{1,0}=-2$. Putting the values from the first two tables in the third we get
\begin{center}
\begin{tabular}{||c c c||} 
 \hline
$(J,K)$ & $t(J,K)$ & $t(J)$ \\ [0.5ex] 
 \hline\hline
$(3,0)$ & $0$ & $4/3\Delta ^3-2/\Delta ^2+1/\Delta $ \\ 
 \hline
$(2,1)$ & $1$ & \\ 
 \hline
$(2,0)$ & $a_1$ & $2(2a_1-1)(1-1/\Delta )/\Delta $ \\  
 \hline
$(1,2)$ & $2$ & \\
 \hline
$(1,1)$ & $8a_1-2$ & \\
 \hline 
$(1,0)$ & $4a_2+4a_1^2-2a_1$ & $(4a_1^2-4a_1+2)/\Delta $  \\
 \hline
\end{tabular}
\end{center}
and Lemma \ref {ail} is done.

\begin {center}
\section* {Proof of Lemma \ref {3}}
\end {center}
This is Theorem 3 (1) of \cite {lauzhao}, but awkwardly our main terms don't immediately match up - if we compare our main term in $E_x(q,a)$ in Lemma \ref {3} with their main term in (1.2) of \cite {lauzhao} then 
\begin {eqnarray*}
\text {our main term }&=&\sum _{r|q}\frac {c_r(a)}{r}\Big (\log (x/r^2)+2\gamma -1\Big )
\\ \text {their main term }&=&\sum _{r|q\atop {r|a}}\frac {\phi (q/r)}{q/r}\Big (\log x/r^2+2\gamma -1\Big )-2\underbrace {\sum _{r|q\atop {r|a}}\sum _{d|q/r}\frac {\mu (d)\log d}{d}}_{(\star )}.
\end {eqnarray*}
But from the observation that
\[ \sum _{dr|q\atop {r|a}}\frac {\mu (d)}{d}=\sum _{d|q}\frac {c_d(a)}{d}\]
we see that these main terms are really the same: the term $(\star )$ is
\begin {eqnarray*}
\sum _{r|q\atop {r|a}}\sum _{d|q\atop {r|d}}\frac {\mu (d/r)}{d/r}\log (d/r)&=&\sum _{d|q}\log d\sum _{r|d\atop {r|a}}\frac {\mu (d/r)}{d/r}-\sum _{r|q\atop {r|a}}\sum _{d|q/r}\frac {\mu (d)}{d}\log r
\\ &=&\sum _{d|q}\frac {\log d}{d}c_d(a)-\sum _{r|q\atop {r|a}}\log r\frac {\phi (q/r)}{q/r}
\end {eqnarray*}
so
\begin {eqnarray*}
\text {their main term }&=&\sum _{r|q\atop {r|a}}\frac {\phi (q/r)}{q/r}\Big (\log x+2\gamma -1\Big )-2\sum _{d|q}\frac {\log d}{d}c_d(a)
\\ &=&\Big (\log x+2\gamma -1\Big )\sum _{r|q\atop {r|a}}\sum _{d|q/r}\frac {\mu (d)}{d}-2\sum _{d|q}\frac {\log d}{d}c_d(a)
\\ &=&\Big (\log x+2\gamma -1\Big )\sum _{d|q}\frac {c_d(a)}{d}-2\sum _{d|q}\frac {\log d}{d}c_d(a)
\\ &=&\text { our main term.}
\end {eqnarray*}

\begin {center}

\end {center}
\hspace {1mm}
\\
\\
\\
\\
\\  
\\ \emph {Tomos Parry
\\ Bilkent University, Ankara, Turkey
\\ tomos.parry1729@hotmail.co.uk}

\end {document}